\documentclass[11pt]{article}

\usepackage{amssymb}
\usepackage{amsmath}
\usepackage{amsthm}
\usepackage{mathrsfs}
\usepackage{graphicx}
\usepackage[utf8]{inputenc}
\usepackage[dvipsnames]{xcolor}

\font\smallit=cmti10
\font\smalltt=cmtt10

\newtheorem{thm}{Theorem}
\newtheorem{prop}{Proposition}
\newtheorem{lem}{Lemma}
\newtheorem{cor}{Corollary}
\newtheorem{rem}{Remark}
\newenvironment{Pf}{\noindent{\bf Proof. }}{\hfill $\blacksquare$ \\}

\def\ts{\textstyle}
\def\lf{\lfloor}
\def\rf{\rfloor}
\def\lc{\lceil}
\def\rc{\rceil}

\def\ov{\overline}

\def\pr{\prime}
\def\sm{\setminus}
\def\niff{\not\leftrightarrow}
\def\a{\alpha}

\def\s{\sigma}
\def\Del{\Delta}
\def\D{\mathscr D}

\begin{document}

\title{\bf On the least size of a graph with a given degree set -- II}
\date{}
\author{
{\bf Jai Moondra}\thanks{{\smallit Department of Computer Science \& Engineering, Indian Institute of Technology, Hauz Khas, New Delhi -- 110016, India} \newline {\smalltt e-mail:\,cs1150230@iitd.ac.in}} 
\qquad 
{\bf Aditya Sahdev}\thanks{{\smallit Department of Computer Science \& Engineering, Indian Institute of Technology, Hauz Khas, New Delhi -- 110016, India} \newline {\smalltt e-mail:\,cs1150207@iitd.ac.in}} 
\qquad
{\bf Amitabha Tripathi}\thanks{{\smallit Department of Mathematics, Indian Institute of Technology, Hauz Khas, New Delhi -- 110016, India} \newline {\smalltt e-mail:\,atripath@maths.iitd.ac.in}}
}
\maketitle

\begin{abstract}
\noindent The {\em degree set\/} of a finite simple graph $G$ is the set of distinct degrees of vertices of $G$. A theorem of Kapoor, Polimeni \& Wall  asserts that the least order of a graph with a given degree set $\D$ is $1+\max \D$. Tripathi \& Vijay considered the analogous problem concerning the least size of graphs with degree set $\D$. We expand on their results, and determine the least size of graphs with degree set $\D$ when (i) $\min \D \mid d$ for each $d \in \D$; (ii) $\min \D=2$; (iii) ${\D}=\{m,m+1,\ldots,n\}$. In addition, given any $\D$, we produce a graph $G$ whose size is within $\min \D$ of the optimal size,  giving a $\big(1+\frac{2}{d_1+1})$-approximation, where $d_1=\max \D$.    
\end{abstract}


\vskip 40pt

\section{Introduction} \label{intro}
\vskip 10pt


A nonincreasing sequence $a_1,\ldots,a_p$ of nonnegative integers is said to be {\em graphic\/} if there exists a simple graph $G$ with vertices $v_1,\ldots,v_p$ such that $v_k$ has degree $a_k$ for each $k$. Any graphic sequence clearly satisfies the conditions $a_k \le p-1$ for each $k$ and $\sum_{k=1}^p\,a_k$ is even. However, these conditions together do not ensure that a sequence will be graphic. Necessary and sufficient conditions for a sequence of nonnegative integers to be graphic are well known; refer \cite{EG60,Hak62,Hav55}. Given a graphic sequence $\mathbf s$ of length $p$, there are polynomial-time algorithms in $p$ to construct a graph with the degree sequence $\mathbf s$; refer \cite{Hak62, Hav55, TVW10}. Therefore, giving a graphic sequence is equivalent to giving some corresponding graph.

The {\em degree set\/} of a simple graph $G$ is the set ${\D}(G)$ consisting of the distinct degrees of vertices in $G$. Conversely, given any set $\D$ of positive integers, a natural question is to investigate the set of all graphs with degree set $\D$, and in particular the least order and size of such graphs.    
The following result answers that question for the least order of a graph with the degree set $\D$. 

\begin{thm}{\bf (\cite{KPW77})} \label{KPW} \\
For each nonempty finite set $\D$ of positive integers, there exists a simple graph $G$ for which ${\D}(G)=\D$. Moreover, there is always such a graph of order $\Del+1$, where $\Del=\max \D$, and there is no such graph of smaller order.
\end{thm} 

\newpage

Tripathi and Vijay \cite{TV06} determine the least size of a graph with a given degree set $\D$, for some special classes of $\D$. We denote this least size by ${\ell}_q(\D)$. In this paper, we determine ${\ell}_q(\D)$ when 
\begin{enumerate}
\item[(i)]
$\min \D \mid d \quad \forall d \in \D$ \\[-22pt]
\item[(ii)] 
$\min \D = 2$ \\[-22pt]
\item[(iii)] 
$\D=\{m,m+1,\ldots,n\}$. \\[-22pt] 
\end{enumerate}
Given any set $\D$, we also give a graph whose size is within ($\min \D-1$) of ${\ell}_q(\D)$, which is a $\big(1+\frac{2}{d_1+1} \big)$-approximation, where $d_1=\max \D$.

Throughout this paper, we let $\D=\{d_1,\ldots,d_n\}$ be a set of positive integers arranged in decreasing order. We shall employ the notation $(d)_m$ to denote $m$ occurrences of the integer $d$. Thus, we may denote a typical sequence by
\begin{equation} \label{notation}
{\mathbf s} = a_1, \ldots, a_p = (d_1)_{m_1}, \ldots, (d_n)_{m_n},
\end{equation}
where $a_j=d_k$ if $\sum_{i = 1}^{k -1} m_i < j \le \sum_{i = 1}^k m_i$, and each $m_k \ge 1$ with $\sum_{i=1}^n m_i = p$. 

We shall write
\[ b_t =\sum_{i=1}^t\:m_i \:\text{for}\: 1 \le t \le n. \]
We call $b_t$ the $t^{\text{th}}$ \emph{breakpoint} of the sequence $\mathbf{s}$, and the set $\{b_t : 1 \le t \le n \}$ as the set of breakpoints of $\mathbf{s}$.

The characterization of graphic sequence due to Erd\H{o}s \& Gallai \cite{EG60} requires verification of as many inequalities as is the order of the graph.

\begin{thm}{\bf (\cite{EG60})} \label{EG_original}
A sequence $\mathbf{s}=a_1,\ldots,a_p$ is graphic if and only if $\:\sum_{i=1}^p a_i$ is even and if the inequalities
\begin{equation*}
\sum_{i = 1}^k a_i \le k(k - 1) + \sum_{i=k+1}^p \min\{a_i, k\}
\end{equation*}
hold for $1 \le k \le p$.
\end{thm}

We also use a refined form of their result, due to Tripathi \& Vijay \cite{TV03} that requires verification of only as many inequalities as the number of distinct terms in the sequence. 

\begin{thm} {\bf (\cite{TV03})} \label{EG} \\
A sequence ${\mathbf s}=a_1,\ldots,a_p$ with the set of breakpoints $\{b_1,\ldots,b_n\}$ is graphic if and only if $\sum_{i=1}^p a_i$ is even and if the inequalities
\[ \sum_{i=1}^k a_i \le k(k-1) + \sum_{i=k+1}^p \min \{k,a_i\} \]
hold for $k \in \{b_1,\ldots,b_n\}$. Moreover, the inequality need only be checked for $1 \le k \le t$, where $t$ is the largest positive integer for which $a_t \ge t-1$.  
\end{thm}
\vskip 5pt

Let sequence ${\mathbf s}$ be as given by eqn.~\eqref{notation}. We set 
\begin{equation} \label{EG_diff}
{\Del}_{\mathbf s}(k) = k(k-1) + \sum_{i=k+1}^p\min \{k,a_i\} - \sum_{i=1}^k a_i, \quad 1 \le k \le p.
\end{equation}
Note that ${\mathbf s}$ is graphic if and only if $\sum_{i=1}^p a_i$ is even and ${\Del}_{\mathbf s}(k) \ge 0$ for $k \in \{b_1,\ldots,b_n\}$ by Theorem~\ref{EG}. 

We denote the sum of the terms of the sequence ${\mathbf s}$ by $\s(\mathbf s)$, and the sum of the elements of the set $S$ by $\s(S)$.   

\section{Basic Results} \label{basic}
\vskip 10pt

In this section, we include two basic results which form the basis of all results in our paper. Given a set $\D$ of positive integers, Proposition \ref{copies} shows the existence of a graphic sequence with exactly one occurrence of each element in $\D$, except that the smallest odd element in $\D$ may occur twice, depending on parity considerations. Additionally, the smallest element in $\D$ will occur multiple times. We also determine the least number of possible occurrences of the smallest element in $\D$ in any such graphic sequence.   

\begin{prop} \label{copies}
Let $\D=\{d_1,\ldots,d_n\}$ be a set of positive integers arranged in decreasing order.

\begin{itemize}
\item[{\rm (a)}] 
Suppose $\s(\D)$ is even or $d_n$ is odd. Let $\mathbf{s}=d_1,\ldots,d_n$, and let
\[ c = \max_{1 \le k \le n-1} \left\{ \left\lc\dfrac{-{\Del}_{\mathbf s}(k)}{\min \{k,d_n\}}\right\rc \right\} = 
           \left\lc\dfrac{-{\Del}_{\mathbf s}(k^{\star})}{\min \{k^{\star},d_n\}}\right\rc. 
\]
Then, there exists a non-negative integer $C$ such that the sequence $\ov{\mathbf{s}}=d_1,\ldots,d_{n-1},(d_n)_{C+1}$ is graphic, and the least such $C$ is given by 
\[ C^{\star} = \begin{cases}
                           c & \mbox{ if $d_n$ and ${\s}(\D)$ are even}; \\ 
                           c & \mbox{ if $d_n$ is odd and ${\s}(\D)+cd_n$ is even}; \\ 
                           c+1 & \mbox{ if $d_n$ and ${\s}(\D)+cd_n$ are odd}.
                          \end{cases} 
\]
Moreover ${\Del}_{\ov{\mathbf s}}(k^{\star})<2d_n$ holds for $C=C^{\star}$.   
\item[{\rm (b)}] 
Suppose $\s(\D)$ is odd and $d_n$ is even. Let $r=\max \{i: d_i \:\text{is odd}\}$, and let ${\mathbf s}=d_1,\ldots,d_{r-1}, (d_r)_2$, $d_{r+1},\ldots, d_n$. Then, there exists a non-negative integer $C$ such that the sequence obtained by appending $C$ copies of $d_n$ to $\mathbf{s}$, $\ov{\mathbf s}=d_1,\ldots, (d_r)_2,d_{r+1},\ldots, (d_n)_{C+1}$ is graphic, and the least such $C$ is given by
\[ C^{\star} = \max_{1 \le k \le n} \left\{ \left\lc\dfrac{-{\Del}_{{\mathbf s}}(k)}{\min \{k,d_n\}}\right\rc \right\} 
    = \left\lc\dfrac{-{\Del}_{\mathbf s}(k^{\star})}{\min \{k^{\star},d_n\}}\right\rc. 
\]
Moreover ${\Del}_{\ov{\mathbf s}}(k^{\star})<d_n$ holds for $C=C^{\star}$. 
\end{itemize}
\end{prop}

\begin{Pf}
\begin{itemize}
\item[{\rm (a)}] 
Suppose ${\s}(\D)$ is even or $d_n$ is odd. Let $k^{\star} \in \{1,\ldots,n-1\}$ be such that $c=\left\lc\dfrac{-{\Del}_{\mathbf s}(k^{\star})}{\min     \{k^{\star},d_n\}}\right\rc$. 
    
Suppose for an arbitrary nonnegative integer $C$, $\ov{{\mathbf s}}=d_1,\ldots,d_{n-1},(d_n)_{C+1}$, so that $b_t=t$ for $1 \le t \le n-1$ and $b_n=n+C$, where $b_t$ is the $t^{\text{th}}$ breakpoint of the sequence $\ov{{\mathbf s}}$. If $C<c$, then 
\begin{eqnarray*}
{\Del}_{\ov{\mathbf s}}(k^{\star}) & = & k^{\star}(k^{\star}-1) + \left( \sum_{i=k^{\star}+1}^n \min \{k^{\star},d_i\} + C \min \{k^{\star},d_n\} \right) - \sum_{i=1}^{k^{\star}} d_i \\
& = & \left( k^{\star}(k^{\star}-1) + \sum_{i=k^{\star}+1}^n \min \{k^{\star},d_i\}  - \sum_{i=1}^{k^{\star}} d_i \right) + C \min \{k^{\star},d_n\} \\
& = &  {\Del}_{\mathbf s}(k^{\star}) + C \min \{k^{\star},d_n\} \\
& \le & {\Del}_{\mathbf s}(k^{\star}) + \left( \left\lc\dfrac{-{\Del}_{\mathbf s}(k^{\star})}{\min \{k^{\star},d_n\}}\right\rc-1 \right) \min                  \{k^{\star},d_n\} \\
& < & {\Del}_{\mathbf s}(k^{\star}) + \left( \dfrac{-{\Del}_{\mathbf s}(k^{\star})}{\min \{k^{\star},d_n\}} \right) \min \{k^{\star},d_n\} \\
& = & 0.
\end{eqnarray*}
Hence $\ov{{\mathbf s}}$ is not graphic when $C<c$ by eqn.~\eqref{EG_diff}.   
    
If $C \ge c$ and $k<n$, then 
\begin{eqnarray*}
{\Del}_{\ov{\mathbf s}}(k) & = & k(k-1) + \left( \sum_{i=k+1}^n \min \{k,d_i\} + C \min \{k,d_n\} \right) - \sum_{i=1}^k d_i \\
& = & \left( k(k-1) + \sum_{i=k+1}^n \min \{k,d_i\}  - \sum_{i=1}^k d_i \right) + C \min \{k,d_n\} \\
& = & {\Del}_{\mathbf s}(k) + C \min \{k,d_n\} \\
& \ge & {\Del}_{\mathbf s}(k) + c \min \{k,d_n\} \\
& \ge & {\Del}_{\mathbf s}(k)  + \left\lc\dfrac{-{\Del}_{\mathbf s}(k)}{\min \{k,d_n\}}\right\rc \min \{k,d_n\} \\
& \ge & {\Del}_{\mathbf s}(k)  + \left( \dfrac{-{\Del}_{\mathbf s}(k)}{\min \{k,d_n\}} \right) \min \{k,d_n\} \\
& = & 0.
\end{eqnarray*}
    
From the definition of $c$ and eqn.~\eqref{EG_diff}, we have 
\[ c \ge -{\Del}_{\mathbf s}(1) = d_1 - (n-1). \]
If $C \ge c$, then $n+C \ge d_1+1$. Thus 
\[ {\Del}_{\ov{\mathbf s}}(n + C) = (n+C)(n+C-1) - \left( \sum_{i=1}^{n-1} d_i + (C+1)d_n \right) \ge (n+C)d_1 - \left( \sum_{i=1}^{n - 1} d_i + (C + 1)d_n \right) \ge 0. 
\] 
Therefore, $\ov{\mathbf s}$ is graphic provided $\s(\ov{\mathbf{s}}) = \big(\sum_{i=1}^{n - 1} d_i \big) + (C + 1)d_n$ is even whenever $C \ge c$. Observe that if $C=c$ and either (i) ${\s}({\mathbf s})$ and $d_n$ are both even, or (ii) $d_n$ is odd and ${\s}({\mathbf s}) + c d_n$ is even, then $\s(\ov{\mathbf{s}})$ is even. Hence $C^{\star}=c$ in these cases. Otherwise, when $d_n$ and ${\s}({\mathbf s})+cd_n$ are both odd, then $\s(\ov{\mathbf{s}})$ is odd when $C=c$ and even when $C=c+1$. Hence $C^{\star}=c+1$ in this remaining case. 
   
Finally, for $C=C^{\star}$ we have 
\begin{eqnarray*}
{\Del}_{\ov{\mathbf s}}(k^{\star}) & = & k^{\star}(k^{\star}-1) + \left( \sum_{i=k^{\star}+1}^n \min \{k^{\star},d_i\} + C^{\star} \min 
\{k^{\star},d_n\} \right) - \sum_{i=1}^{k^{\star}} d_i \\
& \le & \left( k^{\star}(k^{\star}-1) + \sum_{i=k^{\star}+1}^n \min \{k^{\star},d_i\}  - \sum_{i=1}^{k^{\star}} d_i \right) + (c+1) \min 
\{k^{\star},d_n\} \\
& < &  {\Del}_{\mathbf s}(k^{\star}) +  \left( -\dfrac{{\Del}_{\mathbf s}(k^{\star})}{\min \{k^{\star},d_n\}} + 2 \right) \min \{k^{\star},d_n\} \\
& =  & 2 \min \{k^{\star},d_n\} \\
& \le & 2d_n.  
\end{eqnarray*} 
    
\item[{\rm (b)}] 
Notice that $\s(\D)$ is odd implies that $\D$ contains at least one odd integer. Hence $r$ is well defined. Notice also that $\s(\ov{\mathbf s})$ is even for all $C>0$ in this case. Therefore, by Theorem \ref{EG_original} it is enough to show that $\ov{\mathbf{s}}$ satisfies $\Del_{\ov{\mathbf{s}}}(k) \ge 0$, $k \in \{1, \ldots, n - 1, n + C\}$ for $C = C^\star$. The rest of the argument follows along the lines of part (a), and is omitted.
\end{itemize}
\end{Pf}

Given a graph $G$, and a vertex $v$ in $G$, the Splitting lemma allows the construction of a graph $G^{\pr}$ in which the vertex $v$ is replaced by several vertices the sum of degrees of which equals the degree of~$v$. 

\begin{lem} {\bf (Splitting lemma)} \label{splitting} \\[2pt]
Let $G$ be a graph, and let $v \in V(G)$ with $\deg v=d$. Let $(n_1,\ldots,n_r)$ be a partition of $d$ into positive summands. Then there exists a graph $G^{\pr}$ with $V(G^{\pr})=\big(V(G) \sm \{v\}\big) \cup \{v_1,\ldots,v_r\}$ such that $\deg v_i=n_i$, $1 \le i \le r$ and $|E(G^{\pr})|=|E(G)|$. 
\end{lem}

\begin{Pf}
Partition the $d$ neighbours of $v$ into sets $S_1,\ldots,S_r$ with $|S_i|=n_i$, $1 \le i \le r$. Form a graph $G^{\pr}$ from $G$ by replacing the vertex $v$ by vertices $v_1,\ldots,v_r$ such that each $v_i$ is adjacent to the vertices of $S_i$. Note that $|E(G^{\pr})|=|E(G)|$.  
\end{Pf}


\section{Two special cases} \label{divides}
\vskip 10pt

In this section we use the two basic results of Section \ref{basic} to obtain ${\ell}_q(\D)$ when $\min \D$ divides each element of $\D$, and in particular, when $\min \D=1$. We also determine ${\ell}_q(\D)$ when $\min \D=2$. We note that a similar argument may be used to determine ${\ell}_q(\D)$ for $\min \D=3$, but we omit the details.  


\begin{thm} \label{l_q}
Let $\D=\{d_1, \ldots, d_n\}$ be a set of positive integers arranged in decreasing order such that $d_n$ divides $d$ for each $d \in \D$. Then 
\[ {\ell}_q(\D) = \frac{1}{2}\Big( \s(\D) + C^{\star} \big(\min \D\big) \Big), \]
where $\s(\D)$ is the sum of the elements in $\D$ and $C^{\star}$ is as defined in Proposition \ref{copies} (a). 
\end{thm}


\begin{Pf}
Let ${\mathbf s}=(d_1)_{m_1},(d_2)_{m_2},\ldots,(d_n)_{m_n}$ be any graphic sequence with degree set $\D$, so that each $m_i \ge 1$. By repeated applications of Splitting lemma, we may replace all but one copies of $d_i$, $1 \le i<n$, by an appropriate number of $d_n$'s (since $d_n \mid d_i$ for each $i$). Hence we arrive at a graphic sequence $\ov{\mathbf s}=d_1,\ldots,d_{n-1},(d_n)_{M_n}$ for some positive integer $M_n$ such that $\s(\mathbf s)=\s(\ov{\mathbf s})$. Therefore, there exists at least one graphic sequence of the type $\ov{\mathbf s}=d_1,d_2,\ldots,d_{n-1},(d_n)_{M_n}$ for which $\s(\ov{\mathbf s})=2\,{\ell}_q(\D)$. The minimum value of $M_n$ is equal to $C^{\star}+1$, as determined in Proposition \ref{copies} (a).  
\end{Pf}


\begin{cor} \label{min_1}
Let $\D=\{d_1, \ldots, d_n\}$ be a set of positive integers arranged in decreasing order such that $d_n=1$. Then 
\[ {\ell}_q(\D) = \frac{1}{2}\Big( \s(\D) + C^{\star} \Big), \]
where $\s(\D)$ is the sum of the elements in $\D$, and $C^{\star}=\max_{1 \le k \le n} \left\{ -{\Del}_{\mathbf s}(k) \right\}$, where ${\mathbf s}$  is as defined in Proposition \ref{copies} (a). 
\end{cor}

\begin{Pf}
This follows directly from Theorem \ref{l_q}. 
\end{Pf}

\begin{thm} \label{min_2}
Let $\D$ be a set of positive integers such that $\min \D=2$. Then 
\[ {\ell}_q(\D) = \frac{1}{2}\Big( \s(\D) + 2C^{\star} \Big) + 
                               \begin{cases} 
                                  0 & \mbox{ if $\s(\D)$ is even}; \\ 
	                            \frac{1}{2}d_r & \mbox{ if $\s(\D)$ is odd}, 
                               \end{cases}
\]
where $\s(\D)$ is the sum of the elements in $\D$, $r=\max \{i: d_i \:\text{is odd}\}$, and $C^{\star}$ is as defined in Proposition \ref{copies}. 
\end{thm}


\begin{Pf}
This follows directly from Theorem \ref{l_q} when each $d_i$ is even. Henceforth, we assume that at least one $d_i$ is odd. 

Consider the sequence $\ov{\textbf s}$ corresponding to $C=C^{\star}$, as defined in Proposition \ref{copies}. It is easily verified that the sum of the elements of $\ov{\textbf s}$ is given by the expressions for ${\ell}_q(\D)$ in the appropriate cases. 

We prove that any graphic sequence with least element $2$ and degree set $\D$ has size which is at least $\s(\ov{\textbf s})$. Let $\textbf s=(d_1)_{m_1},\ldots,(d_{n-1})_{m_{n-1}},(2)_{m_n}$ be any graphic sequence with degree set $\D$. We may repeatedly apply the Splitting lemma to replace all but one copies of even $d_i$ by $2$'s and all but one copies of odd $d_i$ by one $d_r$ and $2$'s. This gives us a graphic sequence $\textbf u$ with least element $2$, degree set $\D$, with single copies of each $d_i$ except $d_r$ and $2$. 

Let $G$ be a graph with degree sequence $\textbf u$. Assume that $G$ has more than two vertices of degree $d_r$, and let $x$ and $y$ be two of those vertices. Let $v_1,\ldots,v_{d_r}$ and $w_1,\ldots,w_{d_r}$ denote the neighbours of $x$ and $y$ respectively. Without loss of generality, assume $v_1 \ne w_1$ since $d_r>1$. Construct the graph $G^{\pr}$ with vertex set $\big(V(G) \sm \{x,y\}\big) \cup \{a_1,\ldots,a_{\a}\} \cup \{b_1,\ldots,b_{\a}\} \cup \{z\}$, where $\a=\frac{1}{2}(d_r-1)$, and edge set created by removing from the edges of $G$ those with endpoints $x$ or $y$, adding the edges $zv_1$ and $zw_1$, the edges $a_iv_{2i}$ and $a_iv_{2i+1}$, the edges $b_iw_{2i}$ and $b_iw_{2i+1}$ for $i \in \{1,\ldots,\frac{1}{2}(d_r-1) \}$. Thus, 
\begin{eqnarray*}
E\big(G^{\pr}\big) & = & \Big( E(G) \sm \big( \{xv_i: 1 \le i \le d_r\} \:\bigcup\: \{yw_i: 1 \le i \le d_r\} \big) \Big) \:\bigcup\: \{zv_1, zw_1\} \\
& & \bigcup\: \{a_iv_{2i}, a_iv_{2i+1}: 1 \le i \le \ts\frac{1}{2}(d_r-1) \} \:\bigcup\: \{b_iw_{2i}, b_iw_{2i+1}: 1 \le i \le \ts \frac{1}{2}(d_r-1) \}. 
\end{eqnarray*}

We may verify that $G^{\pr}$ has degree set $\D$, with single copies of each $d_i$ except of $d_n=2$ and $d_r$, with the number of copies of $d_r$ decreasing by two, and that $G$ and $G^{\pr}$ have the same size. Repeated applications of this process results in a graph $\tilde{G}$ with degree set $\D$, with single copies of each $d_i$ except of $d_n=2$ and $d_r$, with the number of copies of $d_r$ equal to $1$ or $2$, and that $G$ and $\tilde{G}$ have the same size. The number of copies of $d_r$ is determined by the parity of $\s(\D)$. It now follows from Proposition \ref{copies} that $\tilde{G}$, and consequently $G$, has size at least given by the expressions for ${\ell}_q(\D)$ in the appropriate cases.  
\end{Pf}

\section{The case where $\D$ is an interval} \label{interval}
\vskip 10pt

By an interval we mean a set of the type $\{m, m+1,\ldots, n\}$, where $m,n$ are positive integers with $m \le n$. We denote this by $[m,n]$. In this section, we first determine ${\ell}_q([1,n])$ using Corollary \ref{min_1}, and then use this to determine ${\ell}_q([m,n])$ when $m(m+1)<2\lc \frac{n}{2}\rc$. The case where $m(m+1) \ge 2\lc \frac{n}{2}\rc$ is handled separately.       

\begin{thm} \label{[1,n]}
For $n \ge 1$, 
\[ {\ell}_q\big([1,n]\big) = \frac{1}{2} \left( \frac{1}{2}n(n+1) + \left\lc \frac{n}{2} \right\rc \right). \]
\end{thm}

\begin{Pf}
For the sequence ${\mathbf s}=n,n-1,\ldots,1$ and $k \in \{1,\ldots,n\}$, we may verify that 
\[ {\Del}_{\mathbf s}(k) = \begin{cases} 
                                                -k & \mbox{ if $k<\lceil \frac{n}{2} \rceil$}; \\
                                                2k^2-2(n+1)k+\frac{1}{2}n(n+1) & \mbox{ if $k \ge\lceil \frac{n}{2} \rceil$}. 
                                              \end{cases}
\] 
It is easy to verify that $\max_{1 \le k \le n} \left\{ -{\Del}_{\mathbf s}(k) \right\}=\lceil \frac{n}{2} \rceil$, so the result follows from Corollary  \ref{min_1}. 
\end{Pf}

\begin{prop} \label{m>1}
For any $m \ge 1$, ${\ell}_q\big([m,n]\big) \ge {\ell}_q\big([1,n]\big)$.   
\end{prop}

\begin{Pf}
Let $G$ be any graph of size ${\ell}_q\big([m,n]\big)$ with degree set $[m,n]$. We use the Splitting lemma to construct a graph $G^{\pr}$ with degree set $[1,n]$, and having the same size as $G$, thereby proving the result. 

Choose vertices $v_m,\ldots,v_n$ in $G$ such that $\deg v_i=i$ for each $i \in \{m,\ldots,n\}$, and let $S=\{v_m,\ldots,v_n\}$. Note that $|V(G)| \ge n+1$, so that $|V(G) \sm S| \ge m$. Choose any $m-1$ vertices $w_1,\ldots,w_{m-1}$ in $V(G) \sm S$, and use the Splitting lemma to bifurcate each $w_i$ into $u_i$ and $v_i$ such that $\deg v_i=i$ and $\deg u_i=\deg w_i - i$. Then the graph $G^{\pr}$ with $V(G^{\pr})=\big(V(G) \sm \{w_1,\ldots,w_{m-1}\} \big) \cup \{u_1,\ldots,u_{m-1}\} \cup \{ v_1,\ldots, v_{m-1}\}$ has the same size as $G$, and has degree set $[1,n]$.    
\end{Pf}


\begin{thm} \label{m_small}
Let $m,n$ be positive integers such that $m(m+1) < 2 \lc \frac{n}{2} \rc$. Then
\[ {\ell}_q\big([m,n]\big) = {\ell}_q\big([1,n]\big) = \frac{1}{2} \left( \frac{1}{2}n(n+1) + \left\lc \frac{n}{2} \right\rc \right). \] 
\end{thm}

\begin{Pf}
We first construct a graph $G$ with $\D(G)=[1,n]$ and such that the size of $G$ equals ${\ell}_q\big([1,n]\big)$. Define $V(G)=X \cup Y \cup W$, where 
\[ X = \big\{x_1, \ldots, x_{\lc \frac{n}{2} \rc -1} \big\}, \quad Y = \big\{y_1, \ldots, y_{\lf \frac{n}{2} \rf} \big\}, \quad W = \big\{w, w^{\star} \big\}. \]
Define $E(G)$ by 
\[ \left\{ x_i y_j: 1 \le i \le \left\lc \frac{n}{2} \right\rc -1, 1 \le j \le i \right\} \bigcup \left\{ y_i y_j: 1 \le i < j \le \left\lf \frac{n}{2} \right\rf \right\} 
    \bigcup \big\{ y z: y \in Y, z \in W \big\} \bigcup E_0(G), 
\]
where $E_0(G)$ is empty when $n$ is even, and $E_0(G)=\{w w^{\star}\}$ when $n$ is odd. 

It is easy to verify that $\deg x_i=i$, $1 \le i \le \lc \frac{n}{2} \rc -1$, $\deg y_i=n-i+1$, $1 \le i \le \lf \frac{n}{2} \rf$, and $\deg w=\deg w^{\star}=\lc \frac{n}{2} \rc$. Thus, $G$ is a graph with degree set $[1,n]$ and size ${\ell}_q\big([1,n]\big)=\frac{1}{2} \left( \frac{1}{2}n(n+1) + \left\lc \frac{n}{2} \right\rc \right)$. 

We now construct a graph $G^{\pr}$ with degree set $[m,n]$ and with size equal to the size of $G$. The only vertices in $G$ that have degree less than $m$ are $x_1,\ldots,x_{m-1}$; in fact, $\deg x_i=i$ for $1 \le i \le m-1$. We add $m-i$ edges to the vertex $x_i$ for $i \in \{1,\ldots,m-1\}$ and remove $1+2+3+\cdots+(m-1)=\frac{1}{2}m(m-1)$ edges from the vertex $w^{\pr}$ without affecting the degrees of vertices in $Y$. This results in $\deg x_i=m$ for $i \in \{1,\ldots,m-1\}$, $\deg w^{\star}=\lc \frac{n}{2} \rc - \frac{1}{2}m(m-1)+1 \in [m,n]$, with degrees of all other vertices in $G^{\pr}$ unchanged.     

Let $X^{\pr}=\left\{x_1,\ldots,x_{m-1}\right\}$ and $Y^{\pr}=\left\{y_{\lf \frac{n}{2} \rf -\frac{m(m-1)}{2}+1},\ldots,y_{\lf \frac{n}{2} \rf}\right\}$. Note that $|X^{\pr}|=m-1$, $|Y^{\pr}|=\frac{1}{2}m(m-1)$, and that $x_i \niff y_j$ for $x_i \in X^{\pr}$ and $y_j \in Y^{\pr}$. The nonadjacency in $G$ between vertices of $X^{\pr}$ and $Y^{\pr}$ is a consequence of $\min j=\lf \frac{n}{2} \rf -\frac{1}{2}m(m-1)+1>m-1=\max i$. Note that the degrees of the vertices in $Y^{\pr}$ occupy the $\frac{1}{2}m(m-1)$ consecutive integers starting with $\lc \frac{n}{2} \rc +1$.  

Partition $Y^{\pr}$ into sets $Y_1^{\pr},\ldots,Y_{m-1}^{\pr}$ such that $|Y_i^{\pr}|=i$, $1 \le i \le m-1$. To construct $G^{\pr}$ from $G$, remove the $\frac{1}{2}m(m-1)$ edges $w^{\star} y_j$, $y_j \in Y^{\pr}$, and join $x_i$ to each vertex in $Y_i^{\pr}$, for $i \in \{1,\ldots,m-1\}$.  
It is clear that $G^{\pr}$ has the same size as $G$, that the degree of vertices in $X^{\pr}$ are all equal to $m$, and that the degrees of vertices in $Y^{\pr}$ are unchanged, and $\deg w^{\star}=\lc \frac{n}{2} \rc - \frac{1}{2}m(m-1)+1 \in [m,n]$. Thus $G^{\pr}$ is a graph with desired properties. 

To complete the proof, we note that ${\ell}_q\big([1,n]\big)$ provides a lower bound for ${\ell}_q\big([m,n]\big)$ by Proposition \ref{m>1}, that $G$ and $G^{\pr}$ have size ${\ell}_q\big([1,n]\big)$, and that $\D(G^{\pr})=[m,n]$.   
\end{Pf}

The sequence among 
\begin{equation} \label{sp_seq}
{\mathbf s}_1 = n, n-1, \ldots, m+1, (m)_{m+1}, \qquad {\mathbf s}_2 = n, n-1, \ldots, (m+1)_2, (m)_m 
\end{equation}
with even sum must have optimum size ${\ell}_q([m,n])$ {\em provided} it is graphic. We prove this is the case when $m(m+1) \ge 2 \lc \frac{n}{2} \rc$.  

\begin{thm} \label{m_big}
Let $m,n$ be positive integers such that $m \le n \le 2 \lc \frac{n}{2} \rc \le m(m+1)$. Then exactly one of the sequences ${\mathbf s}_1$, ${\mathbf s}_2$ in eqn.~\eqref{sp_seq} is graphic. In particular, 
\[ {\ell}_q\big([m,n]\big) - {\ell}_q\big([1,n]\big) = \begin{cases} 
                                                                                            \frac{m(m-1)}{4} & \mbox{ if $m \equiv 0, 1\pmod{4}$}; \\ 
                                                                                            \frac{m(m-1)+2}{4} & \mbox{ if $m \equiv 2, 3\pmod{4}$}. 
                                                                                           \end{cases}
\]
\end{thm}

\begin{Pf}
The result is a consequence of \cite[Theorem 4]{TV06} when $m>\lc \frac{n}{2} \rc$. Henceforth, we assume $m \le \lc \frac{n}{2} \rc$.  

Observe that the sum of the integers in ${\mathbf s}_1$ and ${\mathbf s}_2$ differ by one. Hence, {\em at most\/} one of the two sequence can be graphic. We use Theorem \ref{EG} to prove that the sequence with even sum is a graphic sequence. In particular, we show that both sequences satisfy the inequalities in Theorem \ref{EG}. 

Note that 
\[ d_i = \begin{cases} 
                 n-i+1 & \mbox{ if $1 \le i \le n-m$}; \\ 
                 m & \mbox{ if $n-m+1 \le i \le n+1$.} 
               \end{cases}
\]
for sequence ${\mathbf s}_1$, and 
\[ d_i = \begin{cases} 
                 n-i+1 & \mbox{ if $1 \le i \le n-m-1$}; \\
                 m+1 & \mbox{ if $n-m \le i \le n-m+1$}; \\ 
                 m & \mbox{ if $n-m+2 \le i \le n+1$.} 
               \end{cases}
\]
for sequence ${\mathbf s}_2$.

Since $m+1=d_{n-m} \ge n-m-1$ for the sequence ${\mathbf s}_1$ and $m+1=d_{n-m+1} \ge n-m$ for the sequence ${\mathbf s}_2$, while $m=d_{n+1} \le n+1$ for both sequences, the inequality in Theorem \ref{EG} needs to be checked only for $1 \le k \le t$, where $t=n-m$ for the sequence ${\mathbf s}_1$ and $t=n-m+1$ for the sequence ${\mathbf s}_2$. 

Thus, we must show that 
\[ {\Del}_{\mathbf s}(k) = k(k-1) + \sum_{i=k+1}^{n+1}\:\min \{k,d_i\} - \sum_{i=1}^k d_i \ge 0 \;\;\text{for}\;\; 1 \le k \le t \] 
for each of the sequences. 

For the sequence ${\mathbf s}_1$ and $1 \le k \le t$, we may verify that 
\[ {\Del}_{{\mathbf s}_1}(k) = \begin{cases} 
                                                        \frac{1}{2}k(k-1) & \mbox{ if $1 \le k \le m$}; \\
                                                        -k+\frac{1}{2}m(m+1) & \mbox{ if $m < k \le \lceil \frac{n}{2} \rceil$}; \\
                                                        2k^2-2(n+1)k+\frac{1}{2}n(n+1)+\frac{1}{2}m(m+1) & \mbox{ if $\lceil \frac{n}{2} \rceil < k \le t$},
                                                      \end{cases}
\] 
and that ${\Del}_{{\mathbf s}_1}(k) \ge 0$ in each of the three cases. 

The calculations for the sequence ${\mathbf s}_2$ is similar, and omitted here.  
\end{Pf}

\section{An approximation solution for the general case} \label{approx}
\vskip 10pt

In this section, given any set $\D$ of positive integers we construct a graph with degree set $\D$ whose size differs from ${\ell}_q(\D)$ by at most ($\min \D-1$). In particular, we achieve an optimal graph in the case where $\min \D=1$. We show that the sequence constructed in Proposition \ref{copies} achieves this.   


\begin{lem} \label{comparison}
Let $\D=\{d_1,\ldots,d_n\}$ be a set of positive integers arranged in decreasing order. Let $a_1,\ldots,a_p$ denote the graphic sequence $\ov{\mathbf s}$ with degree set $\D$, corresponding to $C=C^{\star}$, as given in Proposition \ref{copies}. If $a_1^{\pr},\ldots,a_{p^{\pr}}^{\pr}$ denotes any graphic sequence with degree set $\D$, then $a_i \le a_i^{\pr}$ for $1 \le i \le \min \{p,p^{\pr}\}$.  
\end{lem}
 
\begin{Pf}
The lemma is clear for the sequence $\ov{\mathbf s}$ in Proposition \ref{copies} (a), that is, when $\s(\D)$ is even or $d_n$ is odd. 


Now assume that $\s(\D)$ is odd and $d_n$ is even, so that the sequence $\ov{\mathbf s}$ is as defined in Proposition \ref{copies} (b). We have $a_i=d_i \le a_i^{\pr}$ for $1 \le i \le r$. Since ${\s}(\D)$ is odd, at least one odd $d_i$, say $d_j$, must be repeated. Since $d_j \ge d_r$, we have $a_i^{\pr} \ge d_{i-1}=a_i$ for $n+1 \ge i \ge r$. The inequality is clear for $i>n+1$. 
\end{Pf}


\begin{thm}\label{approximation-solution}
Let $\D = \{d_1, \ldots, d_n\}$ be a set of positive integers arranged in decreasing order. Then for the graphic sequence $\ov{\mathbf{s}}$ constructed in Proposition \ref{copies} we have
\[ \frac{1}{2} \s\big(\ov{\mathbf{s}}\big) - {\ell}_q(\D) < \min \D. \]
\end{thm}
 
\begin{Pf}
Let the sequence $\ov{\mathbf{s}}=a_1,\ldots,a_p$, as constructed in Proposition \ref{copies}. Let $\ov{\mathbf{s}}^{\pr} = a_1^{\pr}, a_2^{\pr}, \ldots, a_{p^{\pr}}^{\pr}$ be any degree sequence for a graph $G$ with degree set $\D$ with $\ell_q(\D)$ edges. Let $\s\big(\ov{\mathbf{s}}^{\pr}\big) = 2 \ell_q(\D) = \sum_{i=1}^{p^{\pr}}a_i^{\pr}$. 

Thus, we only need to show the following for the result to hold:
\begin{equation} \label{goal}
\s\big( \ov{\mathbf{s}} \big) < \s\big(\ov{\mathbf{s}}^{\pr}\big) + 2d_n. 
\end{equation}

From Lemma \ref{comparison}, we see that
\begin{equation} \label{element_comparison}
a_i \le a_i^{\pr} \qquad \forall \; 1 \le i \le \min\{p, p^{\pr}\}
\end{equation}

Recall that $k^{\star}$ is defined in Proposition \ref{copies} (a), and may be modified suitably for Proposition \ref{copies} (b). We use this notation in either case. Let ${\ell}^{\pr}$ be the largest index such that $a_{{\ell}^{\pr}}^{\pr} \ge k^{\star}$ and $a_{{\ell}^{\pr}+1}^{\pr} < k^{\star}$ and $\ell$ be the largest index such that $a_{\ell} \ge k^{\star}$ and $a_{\ell+1} < k^{\star}$.

It follows from ineq.~\eqref{element_comparison} that
\begin{equation}
\label{l-comparison}
\ell \le {\ell}^{\pr}. 
\end{equation}

Therefore 
\[ \s\big( \ov{\mathbf{s}} \big) = \sum_{i=1}^p a_i = \sum_{i=1}^{k^{\star}} a_i + \sum_{i=k^{\star}+1}^p a_i  = \sum_{i=1}^{k^{\star}} a_i + 
    \sum_{i=k^{\star}+1}^{\ell} \big(a_i - k^{\star} \big) + \sum_{i=k^{\star}+1}^p \min(k^{\star}, a_i).
\]
Let $\s_1=\sum_{i=1}^{k^{\star}} a_i$, $\s_2=\sum_{i=k^{\star}+1}^{\ell} \big(a_i - k^{\star} \big)$, and $\s_3=\sum_{i=k^{\star}+1}^p \min(k^{\star}, a_i)$.

Similarly, 
\[ \s\big( \ov{\mathbf{s}}^{\pr} \big) = \sum_{i=1}^{k^{\star}} a_i^{\pr} + \sum_{i=k^{\star}+1}^{{\ell}^{\pr}} \big(a_i^{\pr} - k^{\star} \big) + 
    \sum_{i=k^{\star}+1}^{p^{\pr}} \min(k^{\star}, a_i^{\pr}).
\]
Let $\s_1^{\pr}=\sum_{i=1}^{k^{\star}} a_i^{\pr}$, $\s_2^{\pr}=\sum_{i=k^{\star}+1}^{{\ell}^{\pr}} \big(a_i^{\pr} - k^{\star} \big)$, and $\s_3^{\pr}=\sum_{i=k^{\star}+1}^{p^{\pr}} \min(k^{\star}, a_i^{\pr})$.
\vskip 5pt

We will show the following three inequalities hold: 
\begin{equation} \label{sigma_ineq}
\s_1 \le \s_1^{\pr}, \quad \s_2 \le \s_2^{\pr}, \quad \s_3 < \s_3^{\pr} + 2d_n
\end{equation}
These imply ineq.~\eqref{goal}. 

From ineq.~\eqref{element_comparison} and ineq.~\eqref{l-comparison}, we get that $\s_1 \le \s_1^{\pr}$ and $\s_2 \le \s_2^{\pr}$. 

To show that $\s_3 \le \s_3^{\pr}$, we apply the Erd\H{o}s-Gallai condition to the graphic sequences $\ov{\mathbf{s}}$ and $\ov{\mathbf{s}}^{\pr}$ at $k^{\star}$, and use Proposition \ref{copies}.  
\[ 0 \le {\Del}_{\ov{\mathbf{s}}}(k^{\star}) = \s_3 - \s_1 - k^{\star}(k^{\star}-1) < 2d_n, \qquad 0 \le {\Del}_{\ov{\mathbf{s}}^{\pr}}(k^{\star}) \le 
    \s_3^{\pr} - \s_1^{\pr} - k^{\star}(k^{\star}-1). 
\]

From ineq.~\eqref{sigma_ineq} and the above inequations, we now have 
\[ \s_3 < \s_1 + k^{\star}(k^{\star}-1) + 2d_n \le \s_1^{\pr} + k^{\star}(k^{\star}-1) + 2d_n \le \s_3^{\pr} + 2d_n. \]

\end{Pf}

\begin{rem}
The sequence $\ov{s}$ in Theorem \ref{approximation-solution} is a $(1+M)$-approximation solution to the problem, where $M=\min \Big\{ \frac{2}{d_1+1}, \frac{2(\sqrt{2}-1)}{n-1} \Big\}$.
\end{rem}

\begin{Pf}
Combining the result of Theorem \ref{approximation-solution} with the basic inequality $2\ell_q(\D)>(d_1+1)d_n$ yields 
\[ \frac{\s\big(\ov{\mathbf{s}}\big)}{2 \ell_q(\D)} < 1 + \frac{d_n}{\ell_q(\D)} < 1 + \frac{2}{d_1+1}. \]
This results in the first bound. 

To obtain the second bound, we use a sharper lower bound $2\ell_q(\D) \ge \left( \sum_{i=1}^n d_i\right)+(d_1+1-n)d_n$ to get
\begin{eqnarray} \label{l_q bound}
2\ell_q(\D) & \ge & \left( \sum_{i=1}^n d_i \right) + (d_1+1-n)d_n \nonumber \\
& \ge & \left( \sum_{i=1}^n d_n + (n-i) \right) + \big( d_1 - (n-1) \big)d_n \nonumber \\
& \ge & nd_n + \frac{1}{2}n(n-1) + d_n^2 \nonumber \\
& \ge & d_n \left( n + \frac{n(n-1)}{2d_n} + d_n  \right) \nonumber \\
& \ge & d_n \left( n + 2\sqrt{\frac{n(n-1)}{2}} \right) \nonumber \\
& > & (1+\sqrt{2}) (n-1)d_n.
\end{eqnarray}
Combining the result of Theorem \ref{approximation-solution} with ineqn.~\eqref{l_q bound} yields 
\[ \frac{\s\big(\ov{\mathbf{s}}\big)}{2 \ell_q(\D)} < 1 + \frac{d_n}{\ell_q(\D)} < 1 + \frac{2(\sqrt{2}-1)}{n-1}. \]
This proves the second bound and hence our claim. 
\end{Pf}
\vskip 5pt

Since our result shows that $2 \ell_q(\D) \in \big(\s(\ov{\mathbf{s}}) - 2\min \D, \s(\ov{\mathbf{s}}) \big]$, a natural algorithm to determine $\ell_q(\D)$ would be to perform a search over this interval. Thus, for each even $\s$ in the interval, we wish to determine if there is a graph $G$ with $\s/2$ edges such that $\D(G)=\D$. One way to do this is to determine determine all solutions to $\s = \sum_{i} m_i d_i$ in positive integers $m_i$, and then check whether any of the corresponding sequences is graphic. This latter problem can be solved in polynomial time in $\sum_i m_i>d_1$ for fixed $m_1, \ldots, m_n$, which is exponential in the input; refer \cite{Hak62, Hav55, TVW10}. The former problem is the well-known Frobenius Coin problem, which has a rich and long history, with several applications and extensions, and connections to several areas of research; refer \cite{Ram05}. The Frobenius Coin problem can be solved in polynomial-time in $d_1 \sum_i m_i$, and is known to be NP-hard \cite{Kan92,Ram05, Ram96}.

Therefore, given a $\s$, we must (i) determine the existence of positive integers $m_1, \ldots, m_n$ such that $\s = \sum_i m_id_i$, and (ii) determine whether $(d_1)_{m_1}, \ldots, (d_n)_{m_n}$ is graphic for each solution $m_1, \ldots, m_n$ in (i). Both these problems have no known polynomial-time algorithms in our input size $\sum_i \log d_i $. While our problem imposes more structure than each of the those problems, we speculate that it is also NP-hard, in which case our result in Theorem \ref{approximation-solution} assumes larger significance.
\vskip 10pt

\noindent {\bf Acknowledgement.} The authors gratefully acknowledge a discussion with Prof. Naveen Garg that led to the results in Section \ref{approx}.

\end{document}